\pdfoutput=1
\RequirePackage{ifpdf}
\ifpdf % We are running pdfTeX in pdf mode   
\documentclass[pdftex]{sigma}
\else
\documentclass{sigma}
\fi

\usepackage{rotating}

\numberwithin{theorem}{section}
\numberwithin{definition}{section}

\begin{document}

\newcommand{\arXivNumber}{1310.7150}

\renewcommand{\vec}[1]{\mathbf{#1}}

\allowdisplaybreaks

\renewcommand{\thefootnote}{$\star$}

\renewcommand{\PaperNumber}{061}

\FirstPageHeading

\ShortArticleName{Twistor Topology of the Fermat Cubic}

\ArticleName{Twistor Topology of the Fermat Cubic\footnote{This paper is a~contribution to the Special Issue on Progress
in Twistor Theory.
The full collection is available at
\href{http://www.emis.de/journals/SIGMA/twistors.html}{http://www.emis.de/journals/SIGMA/twistors.html}}}

\Author{John ARMSTRONG and Simon SALAMON}
\AuthorNameForHeading{J.~Armstrong and S.~Salamon}

\Address{Department of Mathematics, King's College London, Strand, London WC2R 2LS, UK}

\Email{\href{mailto:john.1.armstrong@kcl.ac.uk}{john.1.armstrong@kcl.ac.uk},
\href{mailto:simon.salamon@kcl.ac.uk}{simon.salamon@kcl.ac.uk}}
\URLaddress{\url{https://sites.google.com/site/johnarmstrongmaths/},\\
\hspace*{10.5mm}\url{http://www.mth.kcl.ac.uk/~salamon/}}   

\ArticleDates{Received November 06, 2013, in f\/inal form May 26, 2014; Published online June 06, 2014}

\Abstract{We describe topologically the discriminant locus of a~smooth cubic surface in the complex projective space
${\mathbb{CP}}^3$ that contains $5$~f\/ibres of the projection ${\mathbb{CP}}^3 \longrightarrow S^4$.}


\Classification{53C28; 14N10; 57M20}

\renewcommand{\thefootnote}{\arabic{footnote}}
\setcounter{footnote}{0}

\section{Introduction}

Given an algebraic surface~$\Sigma$ of degree~$d$ in ${\mathbb{CP}}^3$, one can ask how it is conf\/igured relative to the
twistor f\/ibration $\pi:{\mathbb{CP}}^3 \longrightarrow S^4$.
In particular, what can be said about the topology of the restricted map $\pi_{\restriction \Sigma}$?

The motivation for studying this map is two fold.
\begin{enumerate}\itemsep=0pt
\item By the twistor construction, complex hypersurfaces in ${\mathbb{CP}}^3$ transverse to the f\/ibres of~$\pi$
correspond to orthogonal complex structures on $S^4={\mathbb R}^4 \cup \{\infty \}$.
Thus, by studying the interaction of~$\Sigma$ and~$\pi$, we hope to better understand orthogonal complex structures on~${\mathbb R}^4$.

\item The orientation-preserving conformal symmetries of~$S^4$ induce f\/ibre preserving biholomorphisms of~${\mathbb{CP}}^3$.
Thus one can attempt to classify algebraic surfaces in ${\mathbb{CP}}^3$ up to equivalence under the conformal group of~$S^4$ rather than the full group of biholomorphisms of~${\mathbb{CP}}^3$.
Topological invariants of~$\pi_{\restriction \Sigma}$ will automatically give invariants of~$\Sigma$ up to conformal
transformation.
\end{enumerate}

The f\/irst observation that one makes about $\pi_{\restriction \Sigma}$ is that it gives a~$d$-fold branched cover of~$S^4$.
To see this note that each twistor f\/ibre above a~point~$x$ is a~projective line in ${\mathbb{CP}}^3$.
When the def\/ining equation for~$\Sigma$ is restricted to such a~line and written out in homogenous coordinates it either
vanishes or is a~homogenous polynomial $\tilde{p}_x$ in two variables of degree~$d$.
Equivalently we can write it in inhomogeneous coordinates to get a~polynomial $p_x$ of degree $\leqslant d$.
Away from points where $\tilde{p}_x$ vanishes, we have a~$d$-branched cover by the fundamental theorem of algebra.
This motivates the following def\/initions:

\begin{definition}
A {\em twistor line} of~$\Sigma$ is a~f\/ibre of~$\pi$ which lies entirely within the surface~$\Sigma$.
\end{definition}
\begin{definition}
The {\em discriminant locus} of~$\Sigma$ is the branch locus of the map $\pi_{\restriction \Sigma}$.
It is given by the points $x \in S^4$ where either the discriminant of $p_x$ vanishes or the degree of $p_x$ is less
than~$d-1$.
\end{definition}

Thus, given a~surface~$\Sigma$ one can ask how many twistor lines it has, and what is the topology of its discriminant
locus.
These questions were considered in~\cite{salamonViaclovsky} in the case of smooth degree~$2$ surfaces.
It was shown that:
\begin{itemize}\itemsep=0pt
\item The number of twistor lines on a~smooth degree~$2$ surface is either~$0$, $1$, $2$ or~$\infty$.
\item The topology of the discriminant locus for a~degree~$2$ surface is determined by the number of twistor lines.
It is either a~torus ($0$~twistor lines), a~torus pinched at one point ($1$~twistor line), a~torus pinched at two points
($2$ twistor lines) or a~circle ($\infty$ twistor lines).
\end{itemize}
Furthermore, this rough classif\/ication of degree~$2$ surfaces by number of twistor lines was ela\-bo\-rated into a~complete
classif\/ication of degree~$2$ surfaces up to conformal equivalence extending work of Pontecorvo~\cite{pontecorvo}.

In this paper we shall consider the case of cubic surfaces and shall focus in particular upon the Fermat cubic
\begin{gather*}
z_1^3 + z_2^3 + z_3^3 + z_4^3 = 0.
\end{gather*}
Our aim is to compute the topology of the discriminant locus in this particular case.

Note that, as it stands, this problem is not well def\/ined.
The discriminant locus does not depend solely upon the intrinsic geometry of the cubic surface, but also upon the
embedding of the surface in ${\mathbb{CP}}^3$.
Equivalently, it depends upon the choice of twistor f\/ibration.
Fixing the twistor f\/ibration, we will choose an embedding for which the twistor geometry is, in some sense, as simple as
possible.
We choose to measure the ``simplicity'' of an embedding of a~surface~$\Sigma$ f\/irst by the number of twistor lines
on~$\Sigma$ and then by the size of the group of symmetries preserving~$\Sigma$ and the twistor f\/ibration.
In other words, we wish to f\/ind the embedding of the Fermat cubic with the most twistor lines and the largest symmetry
group, and then compute the topology of the discriminant locus in this case.

The remainder of this paper is structured as follows.
In Section~\ref{discriminantLocus} we review the twistor f\/ibration and the calculation of the discriminant locus.
In Section~\ref{twistorLines} we will identify the ``simplest'' embedding of the Fermat cubic relative to the twistor f\/ibration.
We call the resulting embedded cubic surface the ``transformed Fermat cubic''.
Section~\ref{topology} contains the main result of the paper: a~computation of the topology of the discriminant locus of
the transformed Fermat cubic.
This gives the f\/irst known computation of the discriminant locus in degree $d>2$.

This paper is to some extent a~sequel to~\cite{armstrongSalamon} and answers a~question raised in~\cite{povero}.
For background on the twistor f\/ibration we refer readers to~\cite{atiyah}.

\section{Computing the discriminant locus}\label{discriminantLocus}

Let us describe in more detail how one can compute the discriminant locus.
We begin by reviewing an explicit formulation for the twistor projection $\pi:{\mathbb{CP}}^4 \rightarrow S^4$.

Def\/ine two equivalence relations on $({\mathbb H} \times {\mathbb H})\setminus \{0\}$, denoted $\sim_{\mathbb C}$ and
$\sim_{\mathbb H}$, by
\begin{gather*}
(q_1,q_2)\sim_{\mathbb H}(\lambda q_1, \lambda q_2),
\qquad
\lambda \in {\mathbb H} \setminus \{0 \},
\\
(q_1,q_2)\sim_{\mathbb C}(\lambda q_1, \lambda q_2),
\qquad
\lambda \in {\mathbb C} \setminus \{0 \}.
\end{gather*}
We have
\begin{gather*}
(({\mathbb H} \times {\mathbb H}) \setminus \{0\})/{\sim_{\mathbb C}} \cong {\mathbb{CP}}^3
\end{gather*}
and
\begin{gather*}
(({\mathbb H} \times {\mathbb H}) \setminus \{0\})/ {\sim_{\mathbb H}}  \cong {\mathbb HP}^1 \cong {\mathbb R}^4 \cup
\{\infty \} \cong S^4.
\end{gather*}
An explicit map from ${\mathbb{CP}}^3$ to $(({\mathbb H} \times {\mathbb H}) \setminus \{0\})/\sim_{\mathbb C}$ is given
by
\begin{gather*}
[z_1,z_2,z_3,z_4] \rightarrow [z_1+z_2 j, z_3 + z_4j]_{\sim_{\mathbb C}}.
\end{gather*}
The map $\pi:[q_1,q_2]_{\sim_{\mathbb C}} \longrightarrow [q_1,q_2]_{\sim_{\mathbb H}}$ is the twistor f\/ibration.

Left multiplication by the quaternion~$j$ induces an antiholomorphic map $\tau:{\mathbb{CP}}^3 \longrightarrow
{\mathbb{CP}}^3$.
Writing this explicitly we have $\tau[q_1,q_2]_{\sim{\mathbb C}} = [j q_1,j q_2]_{\sim_{\mathbb C}}$.
The map~$\tau$ acts on each f\/ibre of~$\pi$ as the antipodal map. $\tau$ is a~real structure, in other words it is
antiholomorphic and satisf\/ies $\tau^2 = -1$.
Since there is a~unique line through any two points of ${\mathbb{CP}}^3$ one can recover the f\/ibration from~$\tau$
alone: it consists of all projective lines joining a~point~$x$ to the point $\tau x$.

Let us now work out how to compute an explicit expression for the discriminant locus.

Let~$\Sigma$ be a~complex surface in ${\mathbb{CP}}^3$ def\/ined by a~homogeneous polynomial $f_\Sigma(z_1,z_2,z_3,z_4)$
of degree~$d$.
Def\/ine $p_1:{\mathbb R}^4 \rightarrow {\mathbb H} \times {\mathbb H}$ by $p_1(x_1,x_2,x_3,x_4) = [x_1 + x_2 i + x_3 j +
x_4 k, 1]$.
Def\/ine $p_2 = \tau \circ p_1$.
Given a~point $x \in {\mathbb R}^4$, def\/ine $\theta_x:{\mathbb C} \longrightarrow {\mathbb H} \times {\mathbb H} \cong
{\mathbb C} \times {\mathbb C} \times {\mathbb C} \times {\mathbb C}$ by $\theta_x(\lambda) = \lambda p_1(x) + p_2(x)$.
So $\lambda \longrightarrow [\theta_x(\lambda)]_{\sim_{\mathbb C}}$ gives inhomogenous coordinates for the f\/ibre above~$x$.
We see that $f_{\Sigma}(\theta_x(\lambda))$ is a~polynomial $f_x$ in~$\lambda$ with coef\/f\/icients given by polynomials in the~$x_i$.
The discriminant locus is given by the set of points where the degree of $f_x$ is less than $d-1$ or the discriminant of
$f_x$ is equal to~$0$.

Computing this discriminant explicitly, one obtains a~complex-valued polynomial in the coordinates $(x_1,x_2,x_3,x_4)$.
Since the algebraic expression for the discriminant of a~degree~$d$ polynomial always vanishes when evaluated on
a~polynomial of degree ${\leqslant}d-2$, the discriminant locus is given by the zero set of this polynomial.
Taking real and imaginary parts, we obtain an explicit expression for the discriminant locus as the zero set of two real
valued polynomials.

If one runs through this procedure for a~generic cubic surface, the polynomials def\/ining the discriminant locus will be
degree~$12$ polynomials in $4$~variables.
Since the general degree~$12$ polynomial in $4$~variables has $495$~coef\/f\/icients one begins to appreciate the dif\/f\/iculty
in analysing the discriminant locus of a~general cubic surface.
In the next section we will identify the ``simplest'' possible non-singular cubic surface in ${\mathbb{CP}}^3$ relative
to the twistor f\/ibration.

\section{Twistor lines on cubic surfaces}\label{twistorLines}

A beautiful result of classical algebraic geometry due to Salmon states that any smooth cubic surface contains exactly
$27$ projective lines~\cite{cayley}.
When studying the twistor geometry of cubic surfaces it is natural to ask: how many of these lines might be twistor
lines?

The answer to this question can be found by studying the intersection properties of the lines.
The starting point is the classical result, due to Schl\"af\/li~\cite{schlafli} which states that the intersection
properties of the $27$ lines on a~cubic are always the same.
He showed that one can always label six of the lines $a_i$ (with $1\leqslant i \leqslant6$), another six $b_i$ (with
$1\leqslant i \leqslant6$) and the remainder $c_{ij}$ (with $1\leqslant i < j \leqslant6$) in such a~way that:
\begin{itemize}\itemsep=0pt
\item the 6 $a_i$ lines are all disjoint,
\item the 6 $b_i$ lines are all disjoint,
\item $a_i$ intersects $b_j$ if and only if $i \neq j$,
\item $a_i$ intersects $c_{jk}$ if and only if $i \in \{j,k\}$,
\item $b_i$ intersects $c_{jk}$ if and only if $i \in \{j,k\}$,
\item $c_{ij}$ intersects $c_{kl}$ if and only if $\{i,j\} \cup \{j,k\} = \varnothing$.
\end{itemize}
This labelling of the $27$ lines on a~given cubic is far from unique.

It follows immediately that the maximum number of disjoint lines on a~cubic surface is~$6$.
Since twistor lines are necessarily disjoint f\/ibres of~$\pi$ one sees immediately that there can be no more than~$6$
twistor lines on a~cubic surface.
A~more careful analysis allows one to improve upon this.

\begin{theorem}[\cite{armstrongSalamon}]
A~smooth cubic surface contains at most $5$~twistor lines.
If a~smooth cubic surface contains $5$~twistor lines, then the image of these lines under~$\pi$ must all lie on either
a~round $2$-sphere or a~plane in $S^4$.
Note that it follows from Theorem~{\rm \ref{secondTheorem}} below that such cubic surfaces do indeed exist.
\end{theorem}

\begin{proof}
Suppose that~$\Sigma$ is a~cubic surface containing $4$~twistor lines.
Since the lines are skew, it is not dif\/f\/icult to show that one can choose a~labelling of the $27$ lines such that these
lines are~$a_1$, $a_2$, $a_3$, $a_4$.

The line $b_5$ intersects $a_1$, $a_2$, $a_3$, and $a_4$.
Since the twistor lines are invariant under the antiholomorphic involution~$\tau$, the line $\tau(b_5)$ also intersects
$a_1$, $a_2$, $a_3$ and $a_4$.
Thus $\tau(b_5)$ is a~projective line that intersects~$\Sigma$ in $4$~points.
Thus $\tau(b_5)$ is one of the $27$ lines.

Suppose for a~contradiction that $b_5$ intersects $\tau(b_5)$ at a~point~$p$.
Then $b_5$ also intersects $\tau(b_5)$ at $\tau(p)\neq p$.
There is a~unique line through each point so $b_5=\tau(b_5)$.
Hence~$b_5$ is a~f\/ibre of~$\pi$ and so disjoint from~$a_1$.
This is the desired contradiction.

We deduce that $\tau(b_5)$ intersects $a_1$, $a_2$, $a_4$ and $a_4$ and is disjoint from~$b_5$.
By the intersection properties of the cubic, $\tau(b_5)=b_6$.

Suppose that we have a~f\/ifth twistor line~$L$ on~$\Sigma$.
We see that $L \in \{a_5, a_6, c_{56}\}$.
So~$L$ intersects either $b_5$ or $b_6$.
So $L=\tau(L)$ also intersects $\tau(b_5)=b_6$ or $\tau(b_6)=b_5$.
Thus~$L$ inter\-sects~$b_5$ and~$b_6$.
Thus $L=c_{56}$.

Since there is only one possibility for the f\/ifth twistor line, it is impossible to f\/ind a~cubic with six twistor lines.

The lines $b_5$ and $b_6 = \tau(b_5)$ pass through all f\/ive twistor lines.
So~$\pi$ maps the twistor lines onto the image $\pi(b_5)=\pi(b_6)$.
The image of a~projective line under~$\pi$ is always either a~$2$-sphere, two plane or point in~$S^4$.
\end{proof}

Schl\"af\/li used his result on the intersection properties of the $27$ lines as the basis of a~classif\/ication of cubic
surfaces up to projective transformation.
In particular he was able to show that given two skew lines~$b_5$ and~$b_6$, a~generic list of f\/ive points $\{p_i\}$ on~$b_5$
and a~generic list of f\/ive points $\{q_i\}$ on~$b_6$, one can f\/ind a~non-singular cubic surface containing $b_5$
and $b_6$ and the f\/ive lines connecting~$p_i$ to~$q_i$.
Moreover this cubic is unique up to projective transformation.

So given two lines with $b_5=\tau(b_6)$ and $5$~points $p_i$ on $b_5$, one expects to f\/ind a~unique cubic surface up to
projective transformation containing the lines $b_5$, $b_6$ and the lines joining~$p_i$ to~$\tau(p_i)$.
The only issue is that it is not clear that the points $p_i$ and $\tau(p_i)$ will be suf\/f\/iciently generic.
Resolving this question requires a~more detailed study of Schl\"af\/li's results, as was done in~\cite{armstrongSalamon}.
We summarize the main result.

\begin{theorem}\label{secondTheorem}
Let $b_5$ be a~line in ${\mathbb{CP}}^3$ which is not a~twistor fibre.
Let $\{p_1, p_2, p_3, p_4, p_5\}$ be five points on~$b_5$.
We can choose coordinates on~$b_5$ to identify it with the Riemann sphere.

There exists a~smooth cubic surface containing $b_5$, $\tau(b_6)$ and the lines joining $p_i$ and $\tau(p_i)$ if and
only if no four of the points~$p_i$ lie on a~circle under this identification.
\end{theorem}

The data $(b_5, p_1, p_2,
p_3, p_4, p_5)$ in the above theorem uniquely determines the cubic surface up to projective transformation but not up to conformal transformation.
It is shown in~\cite{armstrongSalamon} that there is a~one parameter family of conformally inequivalent cubic surfaces
containing the given lines.
This gives a~classif\/ication of cubic surfaces with f\/ive twistor lines.

So we can now return to the question of determining the ``simplest'' non-singular cubic surface from the point of view
of twistor geometry.
According to the criteria for simplicity given in the introduction, we require the surface to have f\/ive twistor lines
and to have as many conformal symmetries as possible.

Thus we will want to choose f\/ive points on the line $b_5 \cong {\mathbb C} \cup \{\infty \}$ to have as many conformal
symmetries as possible while insisting that no $4$~points lie on a~circle in the Riemann sphere.
Choosing inhomogeneous coordinates for~$b_5$, one sees that the most symmetrical such arrangement is given by the points~$0$,
$1$, $\omega$, $\omega^2$, $\infty$ where~$\omega$ is a~complex cube root of unity.
All equally symmetrical arrangements of f\/ive points on the sphere are conformally equivalent.
The conformal symmetries of these $5$~points are given by the symmetric group~$S_3$: each element of~$S_3$ permutes the
cube roots of unity; odd permutations swap the points~$0$ and~$\infty$.

Up to conformal transformation all choices for the line~$b_5$ are equivalent.
Thus we have narrowed down the search for the non-singular cubic surface with the simplest geometry to a~search among
a~$1$-parameter family of conformally inequivalent surfaces all of which have a~symmetry group of order at least~$6$.
One further symmetry which preserves the twistor f\/ibration is possible: one may have a~symmetry which swaps~$b_5$ to
$\tau(b_6)$ but which leaves the twistor lines f\/ixed.
This allows one to identify the most symmetrical smooth cubic surface with f\/ive twistor lines.
It is unique up to conformal transformation.
An explicit formula for this cubic was found in~\cite{armstrongSalamon}.
It is
\begin{gather}
\label{transformedFermatCubic}
z_1 z_4^2 + z_4 z_1^2 + z_2 z_3^2 + z_3 z_2^2 = 0.
\end{gather}

The simplicity of this formula may lead one to suspect that this is projectively equivalent to a~more familiar cubic
surface.
In fact it is projectively equivalent to the Fermat cubic.
The equivalence is given by the following change of coordinates:
\begin{gather*}
z_1^{\prime}  =  a x_2 + b x_3,
\qquad
z_2^{\prime}  =  c x_2 - b x_3,
\qquad
z_3^{\prime}  =  a x_1 + b x_4,
\qquad
z_4^\prime  =  c x_2 - b x_3,
\end{gather*}
where the constants~$a$,~$b$ and~$c$ are given by $a = \frac{1}{2} + \frac{\sqrt{3}}{6}i$, $b=\overline{a}$ and
$c=\frac{i}{\sqrt{3}}$.

We call the embedded cubic surface def\/ined by equation~\eqref{transformedFermatCubic} the {\em transformed Fermat cubic}
since it is equivalent to the Fermat cubic via a~projective transformation.
Note that this projective transformation is not a~conformal equivalence.
Since the transformed Fermat cubic has $5$~twistor lines and the Fermat cubic only has $3$~these two cubic surfaces are
not conformally equivalent.

\section{The discriminant locus of the transformed Fermat cubic}\label{topology}

Now that we have identif\/ied an appropriately simple cubic surface, we will now compute the topology of its discriminant
locus.

The discriminant locus describes a~singular surface in $S^4 \cong {\mathbb R}^4 \cup \{\infty \}$.
Our strategy to determine its topology is to view it as a~curve in ${\mathbb R}^3$ that changes over time.
We can plot an animation of this curve in Mathematica.
By studying this animation, we hope to deduce the homeomorphism class.
Carefully choosing coordinates for the time axis and for viewing ${\mathbb R}^3$ is essential if one wants to generate
an animation that is easily understood.

Let~$G$ be the group of conformal symmetries of the transformed Fermat cubic.
Viewed in terms of its action on $S_4\cong {\mathbb R}^4 \cup \{\infty\}$,~$G$ is the group of order $12$ generated by
the transformations:
\begin{enumerate}\itemsep=0pt
\item[1)] the rotation~$\theta$ given by
\begin{gather*}
(x_1,x_2,x_3,x_4) \mapsto \big(\cos \big(\tfrac{2\pi}{3}\big) x_1 + \sin\big(\tfrac{2 \pi}{3}\big)x_2, -\sin\big(\tfrac{2\pi}{3}\big) x_1 +
\cos\big(\tfrac{2\pi}{3}\big)x_2, x_3,x_4\big),
\end{gather*}
\item[2)] the
map~$\sigma$ given by $(x_1,x_2,x_3,x_4) \mapsto (x_1, x_2,-x_3, -x_4)$,
\item[3)] the
inversion~$\iota$ given by $ (x_1,x_2,x_3,x_4) \mapsto \tfrac{(x_1,-x_2,-x_3,-x_4)}{x_1^2+x_2^2+x_3^2+x_4^2}$.
\end{enumerate}

{\setlength{\tabcolsep}{0pt}
\begin{sidewaysfigure}[ph!] \centering
\begin{tabular}{cccc}
\includegraphics[scale=0.67]{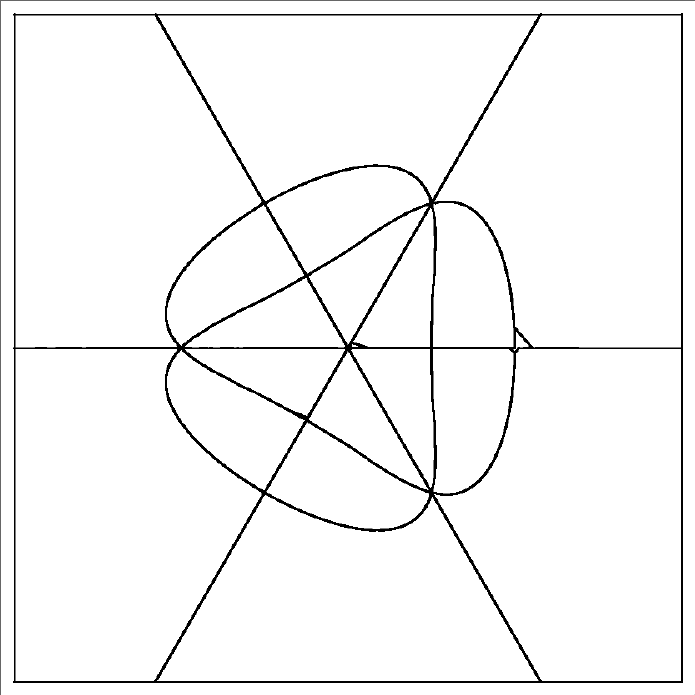} & \includegraphics[scale=0.67]{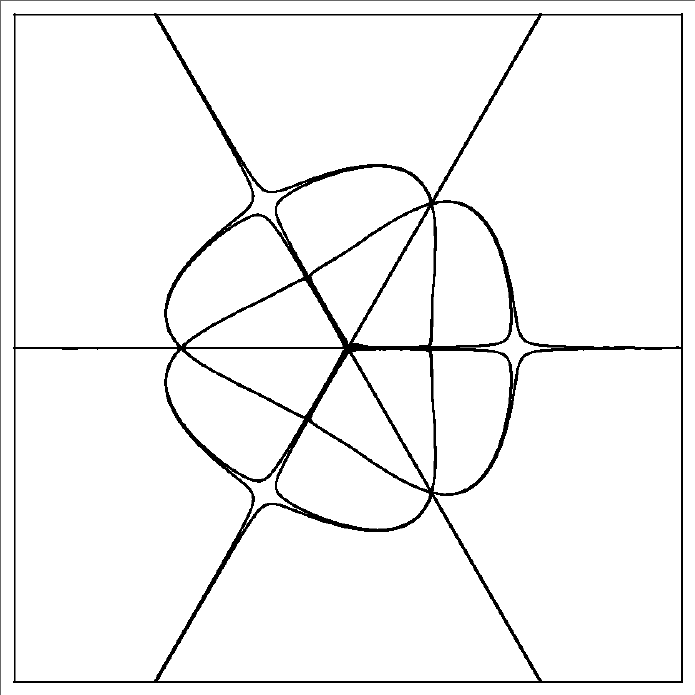} & \includegraphics[scale=0.67]{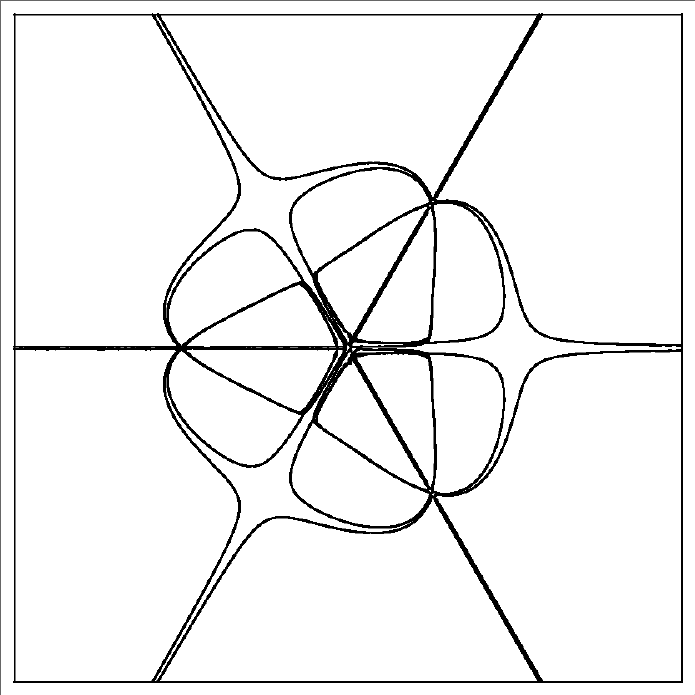} &
\includegraphics[scale=0.67]{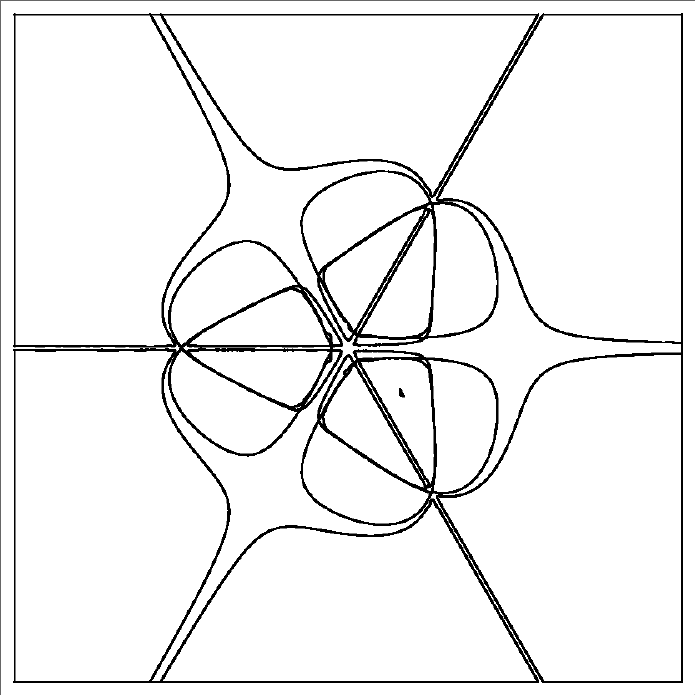}
\\
\includegraphics[scale=0.67]{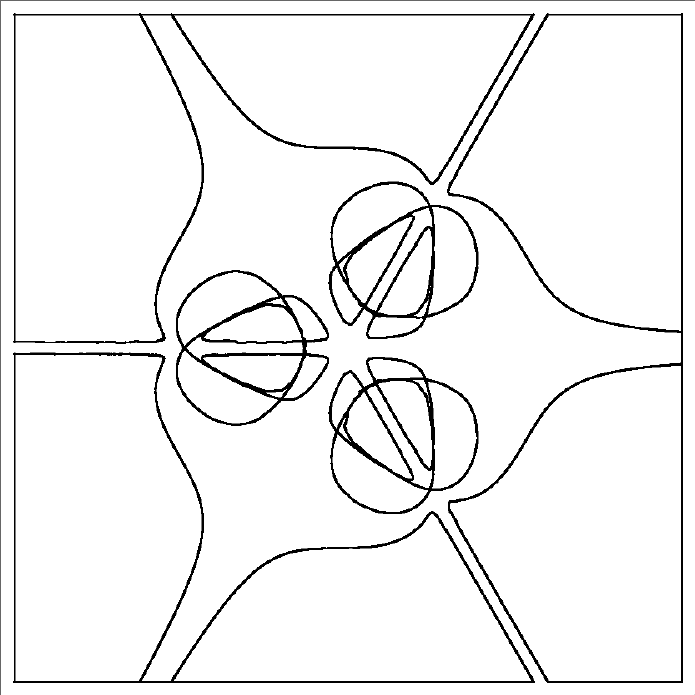} & \includegraphics[scale=0.67]{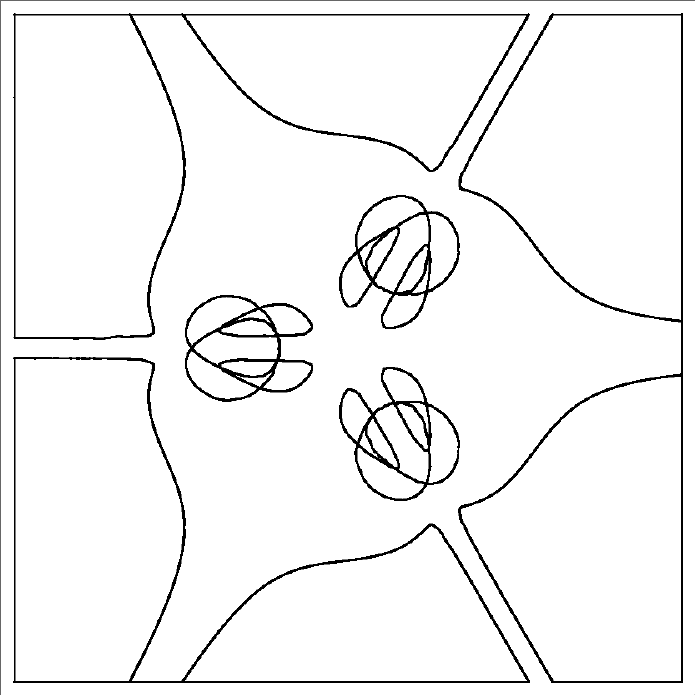} & \includegraphics[scale=0.67]{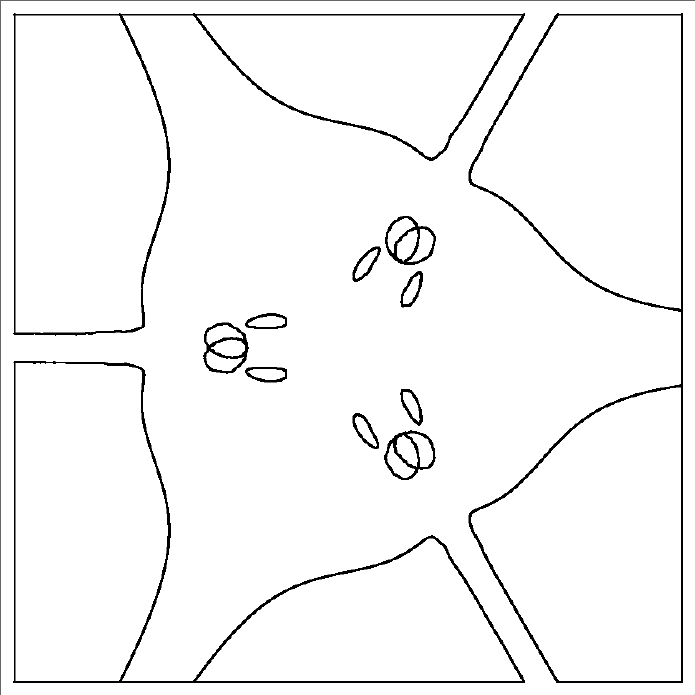} &
\includegraphics[scale=0.67]{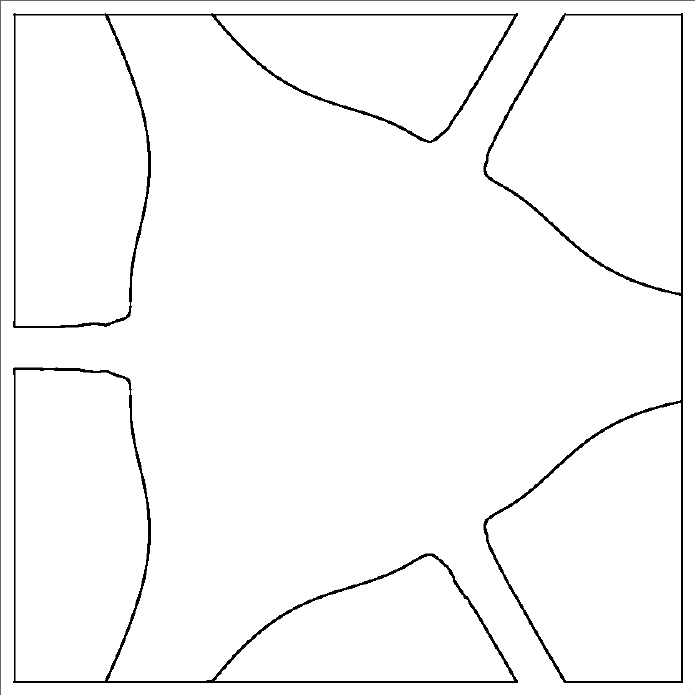}
\end{tabular}

\caption{Animation of the discriminant locus from time $x_4=0$ (top left) to $x_4=0.1$ (bottom right).}
\label{animation}
\end{sidewaysfigure}}

These facts suggest that we use the coordinate $x_4$ as our time coordinate so that the symmetry~$\sigma$ becomes a~time
and parity reversal symmetry.
They also suggest a~perspective from which to view the animation: we should look down the $x_3$ axis onto the $(x_1,x_2)$ plane.

Fig.~\ref{animation} shows
frames of an animation of the discriminant locus from precisely this viewpoint.
It shows the evolution from time $x_4=0$ over positive times.
We have omitted the behaviour over negative times since it can be deduced from the symmetry~$\sigma$.

Fig.~\ref{animation} gives a~rather stylised view of the discriminant locus.
One needs to look at the animation from a~number of angles to understand the situation completely.
However, the main points are already reasonably clear.
At time~$0$, the discriminant locus has a~relatively simple shape consisting of a~number of vertices and curved edges.
At subsequent times, the discriminant locus consists of a~collection of disjoint loops and open curves that leave the
edge of our picture.
The loops do not appear to be disjoint in Fig.~\ref{animation} but that is simply because they overlap when viewed from
this angle.

The loops all shrink to a~point and then disappear (as shown in the f\/ifth and sixth frames).
It is less clear from the picture what will happen to the open curves.
However, we know that the discriminant locus is invariant under the inversion symmetry.
After applying this symmetry, the open curves correspond to disjoint loops.

A loop that shrinks to a~point describes a~surface isomorphic to ${\mathbb R}^2$ in spacetime.
Using the inversion symmetry, we see that the open curves also have worldsheets isomorphic to ${\mathbb R}^2$.
Thus from a~topological perspective, we can think of our animation as providing a~triangulation of the discriminant locus.
The vertices and edges of our triangulation are given by the picture at time~$0$.
The faces are given by the worldsheets of the open curves and loops.

\begin{figure}[ht!] \centering
\includegraphics[scale=0.95]{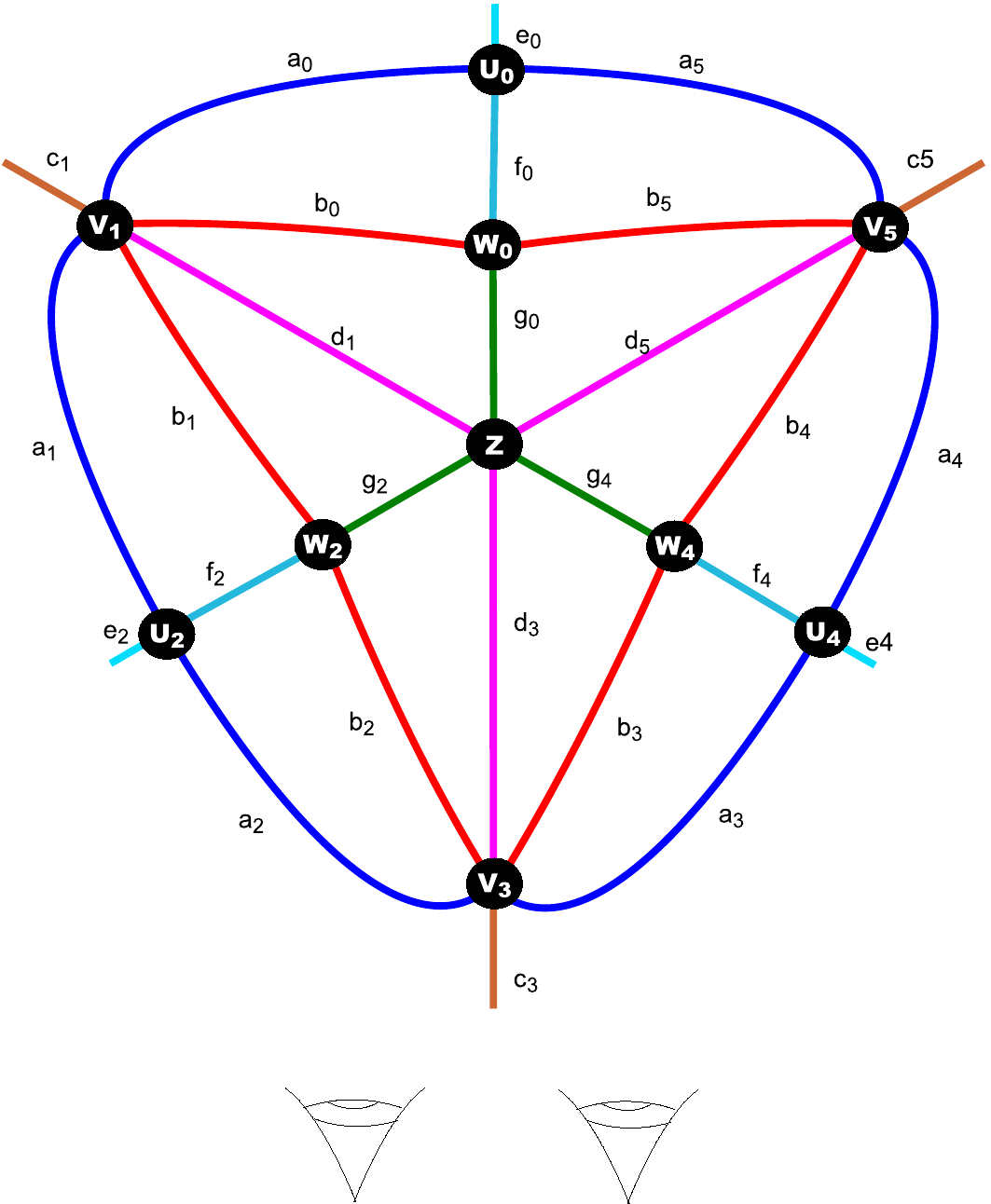}
\caption{The discriminant locus seen from ``above''.}
\label{viewfromabove}
\end{figure}

\begin{sidewaysfigure}[ph!] \centering
\includegraphics{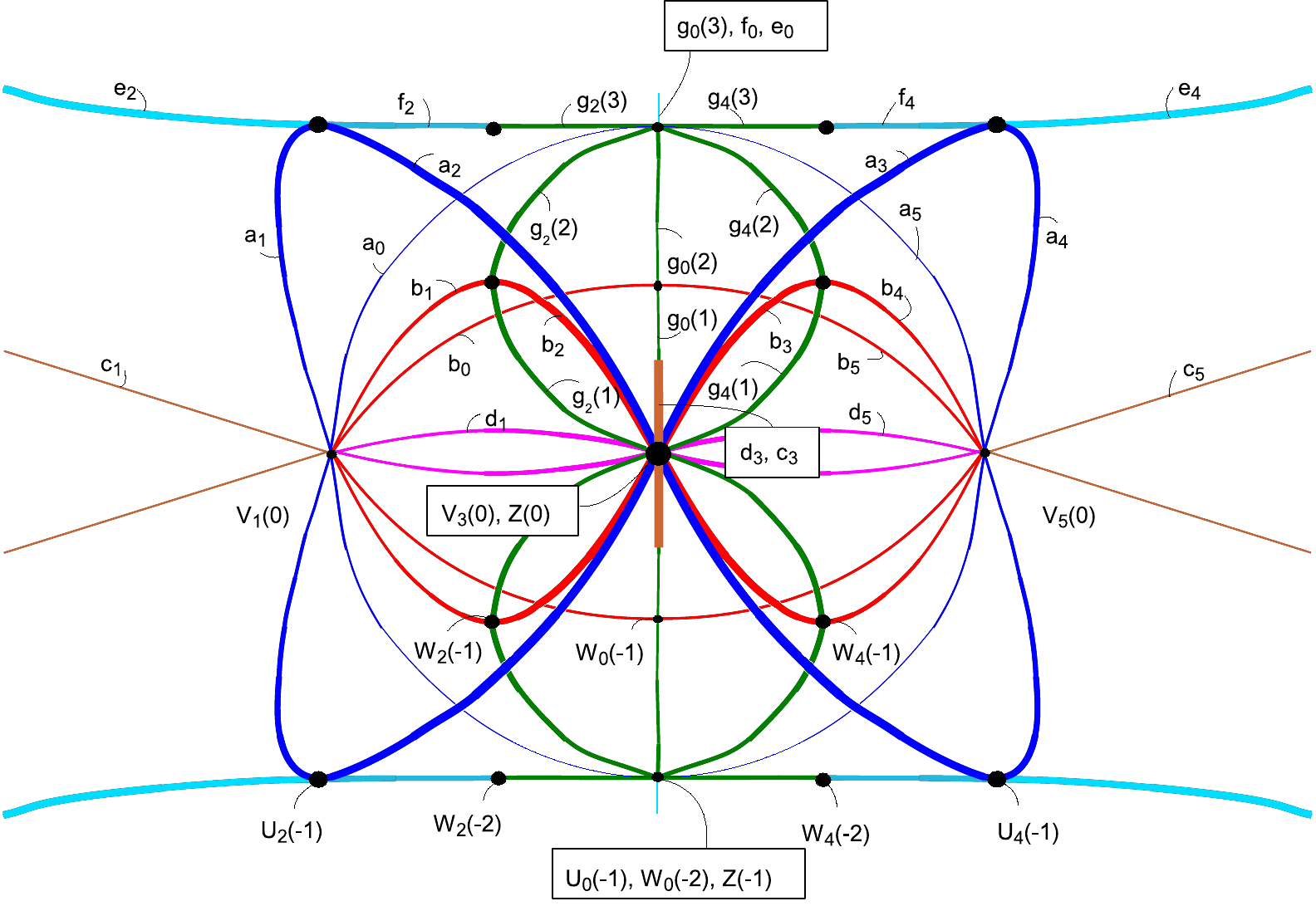}
\caption{The discriminant locus seen from ``the side''.}
\label{viewfromside}
\end{sidewaysfigure}

As we have already remarked, our picture is something of an oversimplif\/ication.
Several curves are overlapping at time~$0$.
We will need to view the discriminant locus from another angle to fully understand what it looks like at time~$0$.
In Figs.~\ref{viewfromabove} and~\ref{viewfromside} we have given two perspectives of the discriminant locus at time
$x_4=0$: a~view from ``above'' showing the $(x_1, x_2)$ plane rotated so that the $x_1$ axis is vertical; a~view from
``the side'' looking at the $(x_2,x_3)$ plane along the $x_1$ axis.
The pair of eyes in Fig.~\ref{viewfromabove} indicates the perspective taken in Fig.~\ref{viewfromside}.

One can identify how the vertices in edges in the view from above correspond to the vertices and edges in the view from
the side.
In Figs.~\ref{viewfromabove} and~\ref{viewfromside} we have labelled all the vertices and edges to show the
correspondence.
There is an obvious action of the group $S_3$ on the view from above which we have incorporated into our labelling.
First we have labelled each $60$ degree sector in the view from above with a~number from~$0$ to~$5$, we have then
labelled each element of the view from above with a~letter that indicates its orbit under the $S_3$ followed by
a~numeric subscript indicating which sector it lies in.
Edges are labelled in lower case, vertices in upper case.
An element such as the edge $g_2$ in the view from above corresponds to a~number of elements in the view from the side
at dif\/ferent heights.
Corresponding to an element in the view from above, we label the elements in the view from the side with a~number in
brackets to indicate the height.
To reduce clutter, in Fig.~\ref{viewfromside} we have only labelled edges in the upper half of the f\/igure and vertices
in the lower half of the f\/igure: the labels in the other half can be deduces by ref\/lecting the picture and changing the
signs in brackets.
For example, the edge $g_2$ in the view from above corresponds to six edges $g_2(3)$, $g_2(2)$, $g_2(1)$, $g_2(-1)$,
$g_2(-2)$ and $g_2(-3)$ in the view from the side.
We have also omitted the number in brackets for elements such as $e_2(1)$ where only the two numbers~$(1)$ and~$(-1)$
would be needed.

One vertex is missing from our picture, a~vertex at inf\/inity.
We label this $Z(\infty)$.
On a~few occasions below we use the notation ${\mathcal E}_m(n)$ for a~generic edge; in such a~formula~$m$ and~$n$ are
integers and $\mathcal E$ is a~symbol like~$a$ or~$b$.

\begin{figure}[ht!]
\centering
\includegraphics[scale=0.53]{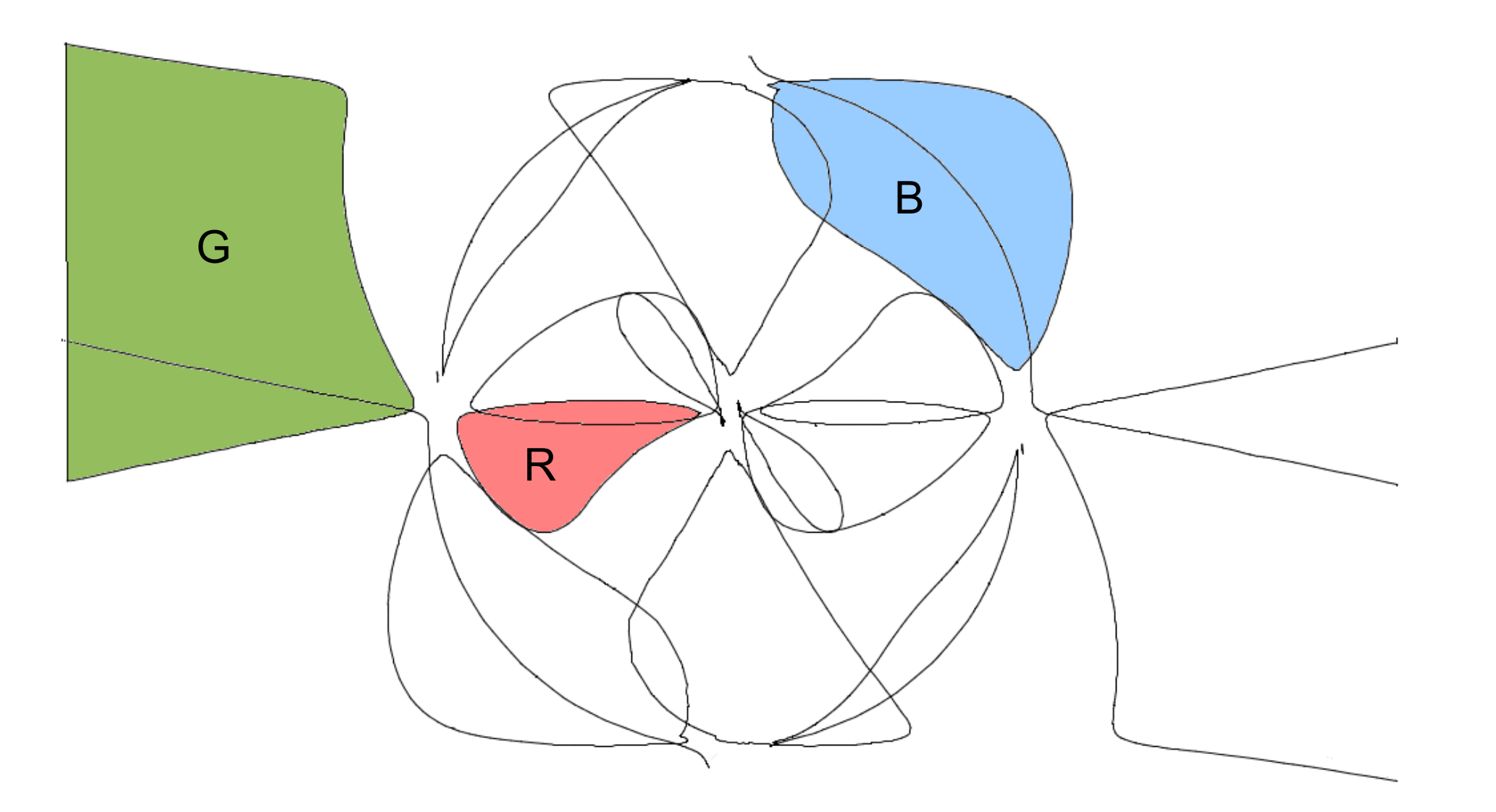} \caption{The view from the side at time $x_4=0.02$.}
\label{timePointNought2}
\end{figure}

We have now identif\/ied the vertices and edges of our triangulation and we have identif\/ied the boundary of each edge.
We now need to identify the faces and the boundary of each face.
To do this we examine the view from the side at time $x_4=0.02$.
This is shown in Fig.~\ref{timePointNought2}.
The value~$0.02$ has been chosen because at this time the curves which correspond to the faces in our triangulation are
very near to the edges of our triangulation.
This makes it easy to read of\/f the boundary corresponding to each face.

Each curve in Fig.~\ref{timePointNought2} corresponds to a~face in our triangulation.
We have labelled three of them as follows:
\begin{enumerate}
\itemsep=0pt
\item[$R$.] The red curve.
The corresponding face in our triangulation has boundary $\{g_2(-1), b_1(-1)$, $d_1(1)\}$.
\item[$G$.] The green curve.
The corresponding face in our triangulation has boundary $\{e_2(1), a_1(1)$, $c_1(-1)\}$.
\item[$B$.] The blue curve.
The corresponding face in our triangulation has boundary $\{f_4(1), g_4(3)$, $g_0(2), b_5(1), -a_4(1))\}$.
\end{enumerate}

\begin{figure}[ht!]
\centering
\includegraphics[scale=1.058]{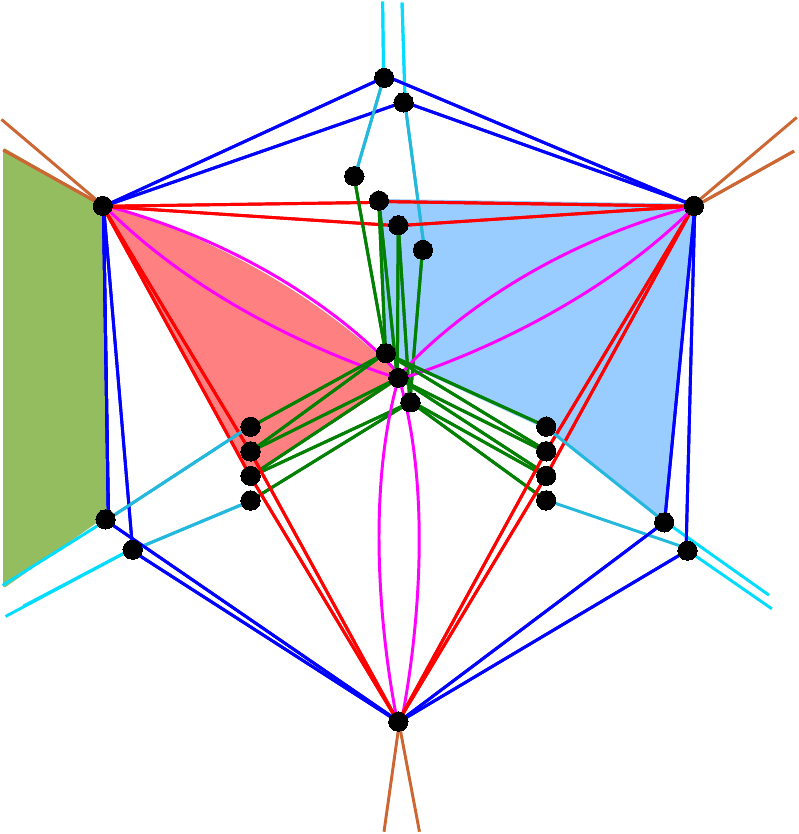} \caption{Schematic view of the triangulation of the discriminant locus.}
\label{schematic}
\end{figure}

We can use the symmetries~$\theta$, $\sigma$ and~$\iota$ to write any face in our triangulation as $\theta^i \sigma^j
\iota^k F$ where $F \in \{R, G, B\}$.
We know that $\theta^i$ sends an edge or vertex ${\cal E}_m(n)$ to ${\cal E}_{(m + 2i \mod 6)}(n)$ and that~$\sigma$
sends ${\cal E}_m(n)$ to ${\cal E}_m(-n)$.
Our notation does not give such a~simple expression for the action of~$\iota$ on edges and vertices.
However, one can check that:
\begin{itemize}\itemsep=0pt
\item $\iota R$ has boundary $\{e_4(1), a_4(1), c_5(-1)\}$,
\item $\iota G$ has boundary $\{g_4(-1), b_4(-1), d_5(1)\}$,
\item $\iota B$ has boundary $\{f_2(-1),g_2(-3), g_0(-2), b_0(-1), a_1(-1)\}$.
\end{itemize}

So one can compute the boundary of the general face $\theta^i \sigma^j \iota^k$ from these formulae.
Thus we have written down all the combinatorial information contained in our triangulation: the list of edges, vertices
and faces and the boundary map $\partial$.
We have found $25$~vertices, $66$~edges and $36$~faces.

The reader may well f\/ind this combinatorial data rather indigestible, but it is summarized schematically in
Fig.~\ref{schematic}.
This shows the view from above perturbed slightly so that all vertices separate.
One should think of this perturbation as corresponding to a~movement along the $x_3$~axis out of the page.
The shaded faces correspond to the~$R$,~$G$ and~$B$ faces shown in the view from the side.

The identif\/ication of this triangulation depends upon a~visual examination of the discriminant locus rather than
a~formal algebraic proof.
If one is willing to accept the results of this visual examination we can now compute the topology of the discriminant
locus.

\begin{figure}[th!]
\centering
\includegraphics[scale=0.9]{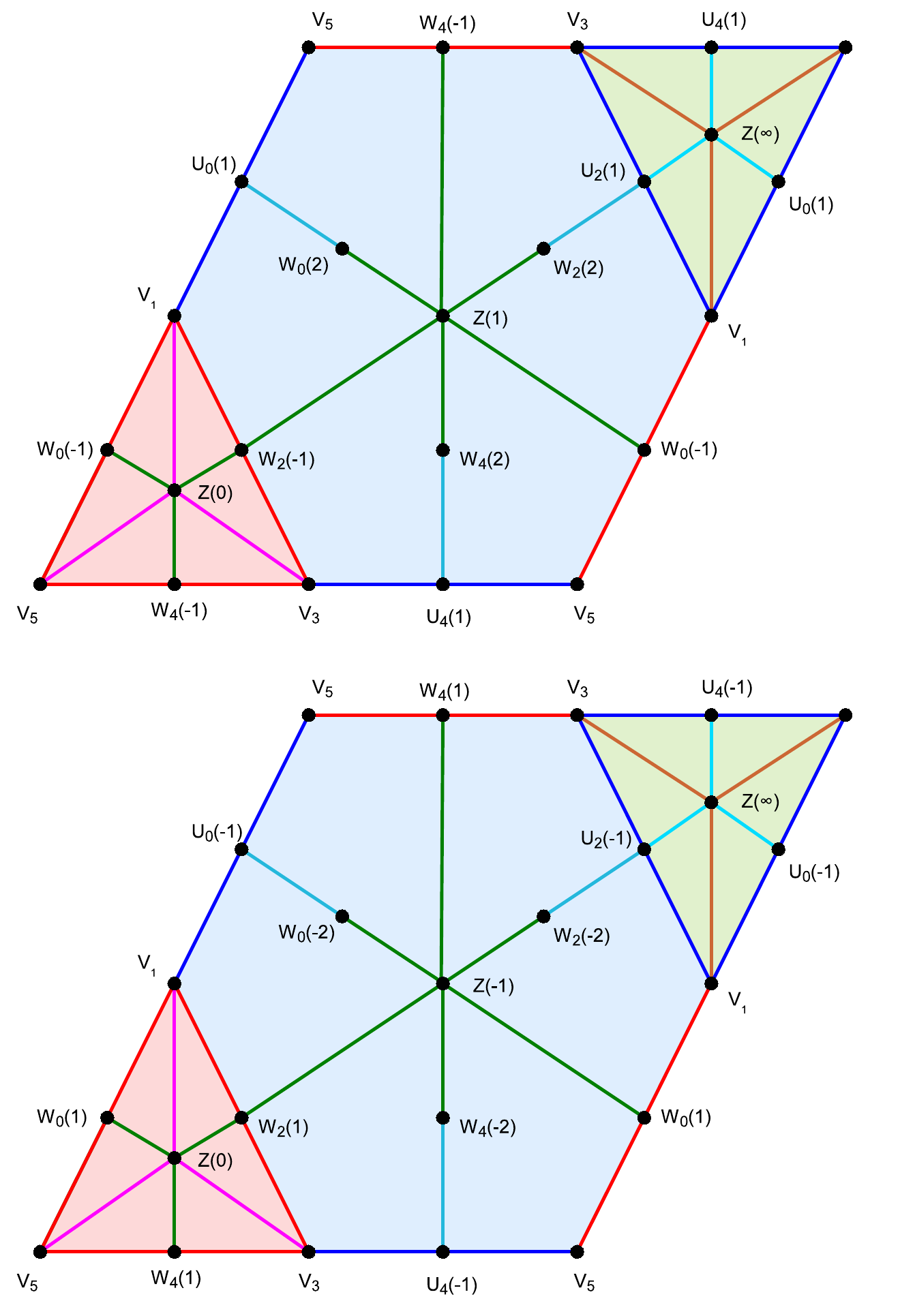} \caption{A triangulation of two tori with $5$~pairs of points identif\/ied.}
\label{twotori}
\end{figure}

\begin{theorem}
The discriminant locus of the transformed Fermat cubic $z_1 z_4^2 + z_4 z_1^2 + z_3 z_2^2 + z_2 z_3^2 = 0$ is
homeomorphic to the space $(T_1 \sqcup T_2)/{\sim}$ where $T_1$ and~$T_2$ are $2$-dimensional tori and~$\sim$ is an
equivalence relation identifying $5$~distinct points on~$T_1$ with another $5$~distinct points on~$T_2$.
\end{theorem}

\begin{proof}
Let~$v$ be an involution def\/ined on the edges of our triangulation by:
\begin{gather*}
v({\cal E}_m(n)) =
\begin{cases}
{\cal E}_{(5-m\, \text{mod}\, 6)}(n), & {\cal E} \in \{a, b \},
\\
{\cal E}_{(-m \,\text{mod}\, 6)}(n), & \hbox{otherwise}.
\end{cases}
\end{gather*}
In terms of Fig.~\ref{schematic} the map~$v$ corresponds to a~ref\/lection of the edges in the vertical, $x_1$, axis.

Let $\partial$ denote the map sending an oriented surface to its boundary.
Notice that $\partial (\iota R) = v \partial G$, $\partial (\iota G) = v \partial R$ and $\partial (\iota B) = \sigma v
\partial B$.
Thus the set of all boundaries of faces in our triangulation is generated by the action of~$\theta$,~$\sigma$ and~$v$ on
the boundaries of~$R$,~$G$ and~$B$.
This means that the triangulation of the discriminant locus is equivalent to another triangulation with edges also given
by the schematic diagram~\ref{schematic} but with faces generated by the faces in~\ref{schematic} under the
maps~$\theta$,~$\sigma$ and~$v$.
The point is that we have eliminated the rather awkward inversion symmetry~$\iota$ and replaced it with the simple
ref\/lection~$v$.

Using these symmetries, it is easy to check that our triangulation is isomorphic to that shown in Fig.~\ref{twotori}
once opposite edges of the parallelograms have been identif\/ied to obtain tori and the vertices
\begin{gather}
Z(0),
\qquad
Z(\infty),
\qquad
V_1,
\qquad
V_3,
\qquad
V_5
\label{identifiedPoints}
\end{gather}
have been identif\/ied.
The mappings of edges and faces under the isomorphism can be easily deduced.

Thus the discriminant locus is homeomorphic to two tori with f\/ive pairs of points identif\/ied.
The f\/ive points~\eqref{identifiedPoints} are the images of the f\/ive twistor lines in the cubic~$\Sigma$.
\end{proof}

The Mathematica code used to generate the f\/igures in this article is available from
\href{https://kclpure.kcl.ac.uk/portal/en/publications/twistor-topology-of-the-fermat-cubic--mathematica-code(0f8cebf5-a9f1-41cb-8511-e79f54a1c819).html}{https:// kclpure.kcl.ac.uk/portal/en/publications/twistor-topology-of-the-fermat-cubic{-}{-}mathematica-\linebreak code(0f8cebf5-a9f1-41cb-8511-e79f54a1c819).html}.

\pdfbookmark[1]{References}{ref}
\LastPageEnding

\end{document}